\documentclass[11pt]{article}
\usepackage{amssymb,amsmath}
\usepackage{amsmath}
\usepackage{amsfonts}
\newtheorem{precor}{{\bf Corollary}}

\newtheorem{precon}{{\bf Conjecture}}

\newtheorem{predefin}{{\bf Definition}}

\newtheorem{preexm}{{\bf Example}}

\newtheorem{preappl}{{\bf Application}}

\newtheorem{prelem}{{\bf Lemma}}

\newtheorem{preproof}{{\bf Proof.\ }}

\newenvironment{proof}[1]{\begin{preproof}{\rm
               #1}\hfill{$\blacksquare$}}{\end{preproof}}
\newtheorem{presproof}{{\bf Sketch of Proof.\ }}

\newtheorem{prethm}{{\bf Theorem}}

\newenvironment{thm}{\begin{prethm}{\hspace{-0.5
               em}{\bf.\ }}}{\end{prethm}}
\newtheorem{prealphthm}{{\bf Theorem}}

\newtheorem{prepro}{{\bf Proposition}}

\newtheorem{preprb}{{\bf Problem}}

\def\conct[#1,#2]{\mbox {${#1} \leftrightarrow {#2}$}}
\def\dconct[#1,#2]{\mbox {${#1} \rightarrow {#2}$}}
\def\deg[#1,#2]{\mbox {$d_{_{#1}}(#2)$}}
\def\mindeg[#1]{\mbox {$\delta_{_{#1}}$}}
\def\maxdeg[#1]{\mbox {$\Delta_{_{#1}}$}}
\def\outdeg[#1,#2]{\mbox {$d_{_{#1}}^{^+}(#2)$}}
\def\minoutdeg[#1]{\mbox {$\delta_{_{#1}}^{^+}$}}
\def\maxoutdeg[#1]{\mbox {$\Delta_{_{#1}}^{^+}$}}
\def\indeg[#1,#2]{\mbox {$d_{_{#1}}^{^-}(#2)$}}
\def\minindeg[#1]{\mbox {$\delta_{_{#1}}^{^-}$}}
\def\maxindeg[#1]{\mbox {$\Delta_{_{#1}}^{^-}$}}

\def\dre[#1,#2,#3]{\mbox {${\cal E}_{_{#3}}(#1,#2)$}}
\def\pdre[#1,#2,#3]{\mbox {${\cal P}_{_{#3}}(#1,#2)$}}
\def\var[#1,#2]{\mbox {${\rm Var}_{_{#1}}(#2)$}}
\def\ls[#1]{\mbox {$\xi^{^{#1}}$}}
\def\hom[#1,#2]{\mbox {${\rm Hom}({#1},{#2})$}}
\def\onvhom[#1,#2]{\mbox {${\rm Hom^{v}}(#1,#2)$}}
\def\onehom[#1,#2]{\mbox {${\rm Hom^{e}}(#1,#2)$}}
\def\core[#1]{\mbox {$#1^{^{\bullet}}$}}
\def\cay[#1,#2]{\mbox {${\rm Cay}({#1},{#2})$}}
\def\cays[#1,#2]{\mbox {${\rm Cay_{s}}({#1},{#2})$}}
\def\dirc[#1]{\mbox {$\stackrel{\rightarrow}{C}_{_{#1}}$}}
\def\cycl[#1]{\mbox {${\bf Z}_{_{#1}}$}}

\def\sdg[#1]{\mbox {$\stackrel{\leftrightarrow}{#1}$}}
\begin{document}
%\setcounter{page}{183}
%{\footnotesize AAA {\bf ?} (200?) ?--?}\\
%\maketitle
\begin{center}
{\Large \bf A Note on Chromatic Sum}\\
\vspace*{0.5cm}
{\bf Meysam Alishahi and Ali Taherkhani\footnote{Corresponding author. Tel.: +98 2129902917.} }\\
{\it Department of Mathematical Sciences}\\
{\it Shahid Beheshti University, G.C.,}\\
{\it P.O. Box {\rm 19839-63113}, Tehran, Iran}\\
{\tt m\_alishahi@sbu.ac.ir}\\
%{\tt hhaji@sbu.ac.ir}\\
{\tt a\_taherkhani@sbu.ac.ir}\\ \ \\
\end{center}
\footnotetext[2]{This paper is partially supported by Shahid
Beheshti University.}
\begin{abstract}
\noindent The chromatic sum $\Sigma(G)$ of a graph $G$ is the
smallest sum of colors among of proper coloring with the natural
number. In this paper, we introduce a necessary condition for the
existence of graph homomorphisms. Also,  we present
$\Sigma(G)<\chi_f(G)|G|$ for every graph $G$.

\begin{itemize}
\item[]{{\footnotesize {\bf Key words:}\ chromatic sum, graph homomorphism, Fractional chromatic number.}}
\item[]{ {\footnotesize {\bf Subject classification: 05C} .}}
\end{itemize}
\end{abstract}
\section{Introduction and Preliminaries}
We consider finite undirected graphs with no loops and multiple
edges and use \cite{MR2089014}  for the notions and notations not
defined here. Let $G$ be a graph and $c$ be a proper coloring of
it, define $\Sigma_c(G)=\sum_{\{v\in V(G)\}}c(v)$. The
vertex-chromatic sum of $G$, denoted by $\Sigma(G)$, is defined
as $\min\{\Sigma_c(G)|\ c\ {\rm is\ a\ proper\ coloring\ of\
}G\}$. The vertex-strength of G denoted by $s(G)$, or briefly by
$s$, is the smallest number $s$ such that there is a proper
coloring $c$ with $s$ colors where $\Sigma_c(G) = \Sigma(G)$.
Clearly, $s(G) \geq \chi(G)$ and equality does not always hold.
In fact, for every positive integer $k$, almost all trees satisfy
$s>k$; see \cite{kobik}. Chromatic sum has been investigated in
literature
\cite{MR1041612,MR1746464,tabular,MR1722229,MR2085077,kobik,MR1000084}.

In \cite{MR1000084}, Thomassen et al. obtained several bounds for
chromatic sum for general graphs. The first is a rather natural
result of an application of a greedy algorithm: $\Sigma(G)\leq
n+e$, where $n$ and $e$ are the number of vertices and edges of
$G$, respectively. Also, they presented an upper and lower limit
for the chromatic sum in terms of $e$. They showed that
$\sqrt{8e}\leq\Sigma(G)\leq\frac32(e+1)$ and these bounds are
sharp.

Let $G$ and $H$ be two graphs. A  homomorphism $\sigma$ from a
graph $G$ to a graph $H$ is a map $\sigma : V(G) \longrightarrow
V(H)$ such that $uv \in E(G)$ implies $\sigma(u)\sigma(v) \in
E(H)$. The set of all homomorphisms from $G$ to $H$ is denoted by
$\hom[G,H]$. An isomorphism of $G$ to $H$ is a homomorphism
$f:G\rightarrow H$ which is a vertex and edge bijective
homomorphism. An isomorphism $f:G\rightarrow G$ is called an
automorphism of $G$, and the set of all automorphism of $G$ is
denoted by ${\rm Aut}(G)$.

Suppose $m \geq 2n$ are positive integers. We denote by $[m]$ the
set $\{1, 2, \cdots, m\}$, and denote by ${[m] \choose n}$ the
collection of all $n$-subsets of $[m]$. The  Kneser graph
$KG(m,n)$ has vertex set ${[m] \choose n}$, in which $A \sim B$
if and only if $A \cap B = \emptyset$. The graph $KG(5,2)$ is
named Petersen graph that is denoted by $P$. It was conjectured by
Kneser in 1955 and proved by Lov\'{a}sz \cite{MR514625} in 1978
that $\chi(KG(m,n))=m-2n+2$.

The fractional chromatic number of a graph $G$, denoted by
$\chi_f(G)$, is the infimum of the ratios $\frac mn$ such that
there is a homohomrphism from $G$ to $KG(m,n)$. It is known
\cite{MR1481157} that the infimum in the definition can be
attained, and hence can be replaced by the minimum. It is easy to
see $\chi_f(G)\leq \chi(G)$. On the other hand, the ratio
$\frac{\chi(G)}{\chi_f(G)}$ can be arbitrary large, see
\cite{MR1481157}.

In next section we present a necessary condition for existence of
graph homomorphisms in terms of chromatic sum. Next, we introduce
an upper bound for chromatic sum based on fractional chromatic
number.

%%%%%%%%%%%%%%%%%%%%%%%%%%%%%%%%%%%%%%%%%%%%%%%%%%%%%%%%%%%%%%%%%%%%%%%%%%%%%%
%%%%%%%%%%%%%%%%%%%%%%%%%%%%%%%%%%%%%%%%%%%%%%%%%%%%%%%%%%%%%%%%%%%%%%%%%%%%%%
\section{Graph Homomorphism and Chromatic Sum}
Graph homomorphism is a fundamental concept in graph theory,
where it is related to many important concepts and problems in
the field. It is well-known that in general it is a hard problem
to decide whether there exists a homomorphism from a given graph
$G$ to a given graph $H$, and consequently, it is interesting to
obtain necessary conditions for the existence of such mappings. In
this regard, we have the following theorem.

\begin{thm}\label{hom}
Let $G$ and $H$ be two graphs such that $H$ is a vertex transitive
graph. If $\sigma:G\longrightarrow H$ is a homomorphism, then
$$\frac{\Sigma(G)}{|G|}\leq\frac{\Sigma(H)}{|H|}.$$
\end{thm}
\begin{proof}{
Let ${\rm Aut(H)}=\{f_1,f_2,\ldots,f_t\}$ and
$\tilde{G}=\bigcup_{i=1}^t G_i$ that $G_i$ is an isomorphic  copy
of $G$. Define $\tilde{\sigma}:\tilde{G}\longrightarrow H$ such
that its restriction to $G_i$ is $f_i\circ\sigma$. Since $H$ is a
vertex transitive graph, one can easily show that for every $v\in
V(H)$, $|\tilde{\sigma}^{-1}(v)|=t\frac{|G|}{|H|}$ and it is
independent of $v$. Now, suppose $c$ is a proper coloring of $H$
such that $\Sigma_c(H)=\Sigma(H)$. For any vertex $v\in V(\tilde
G)$, set $\tilde{c}(v)=c(\tilde{\sigma}(v))$. Obviously,
$\tilde{c}$ is a proper coloring of $\tilde{G}$ and also
$\Sigma_{\tilde{c}}(\tilde G)=\frac{t|G|}{|H|}\times\Sigma(H)$.
Therefore, there is an $i$ such that
$\Sigma_{\tilde{c}_{|G_i}}(G_i)\leq
\frac{|G|}{|H|}\times\Sigma(H)$ and since $G=G_i$,  $\Sigma(G)\leq
\frac{|G|}{|H|}\times\Sigma(H)$ which is the desired conclusion. }
\end{proof}

Theorem \ref{hom} provides a necessary condition for the
existence of graph homomorphisms. Here we show that The Petersen
graph $P$ has the same chromatic number  and circular chromatic
number. One can check that $\Sigma(P)=19$ and
$\Sigma(K_{\frac83})=15$. Therefore, as an application of the
previous theorem, there is no homomorphism from $P$ to
$K_{\frac83}$.

It is well-known that the chromatic sum  is an NP-complete
problem\cite{kobik}. In this regard, finding upper and lower
bounds for chromatic sum is useful. It was shown in \cite{tabular}
that $\Sigma(G)\leq(\frac{\chi(G)+1}{2})|G|$. Since
$\Sigma(K_n)=\frac{n(n+1)}{2}$, if we set $H=K_{\chi(G)}$, then
Theorem \ref{hom} implies  this bound.  Here we obtain an upper
bound for the chromatic sum in terms of fractional chromatic
number.

For an independent set $S$ in a graph $G$ the following
inequality is an  immediate consequence of the definition of the
chromatic sum (\cite{MR1746464}),
\begin{equation}\label{lem}
\Sigma(G)\leq |G|+\Sigma(G\setminus S).
\end{equation}

%%%%%%%%%%%%%%%%%%%%%%%%%%%%%%%%%%%%%%%%%%%%%%%%%%%%%%%%%%%%%%%%%%%%%%%%%%%%%%
%%%%%%%%%%%%%%%%%%%%%%%%%%%%%%%%%%%%%%%%%%%%%%%%%%%%%%%%%%%%%%%%%%%%%%%%%%%%%%
\begin{thm}\label{chif}
For every graph $G$, we have $$\Sigma(G)<\chi_f(G)|G|.$$
\end{thm}
\begin{proof}{
Assume that $\chi_f(G)=\frac mn$ and ${\rm
Hom}(G,KG(m,n))\neq\emptyset$. In view of equation \ref{lem}, we
have $\Sigma(KG(m,n))\leq {m\choose n}+\Sigma(KG(m-1,n))$. Hence
$\Sigma(KG(m,n))\leq \sum_{i=0}^{m-2n-1}{m-i\choose
n}+\Sigma(KG(2n,n)).$ On the other hand,
$\sum_{i=0}^{m-2n-1}{m-i\choose n}={m+1\choose n+1}-{2n+1\choose
n+1}$ and $\Sigma(KG(2n,n))=\frac32{2n\choose n}$. Therefore,
$\Sigma(KG(m,n))\leq {m+1\choose n+1}-(\frac{n-1}{2n+2}){2n
\choose n}.$ Now, since ${\rm Hom}(G,KG(m,n))\neq\emptyset$,
Theorem \ref{hom} implies that
$$
\Sigma(G)\leq\left(\frac{m+1}{n+1}-(\frac{n-1}{2n+2})\frac{{2n
\choose n}}{{m \choose n}}\right)|G|.
$$
Furthermore, $\frac{m+1}{n+1}-(\frac{n-1}{2n+2})\frac{{2n \choose
n}}{{m \choose n}} \leq \frac {m+1}{n+1}<\frac mn=\chi_f(G)$, as
desired. }
\end{proof}
%%%%%%%%%%%%%%%%%%%%%%%%%%%%%%%%%%%%%%%%%%%%%%%%%%%%%%%%%%%%%%%%%%%%%%%%%%%%%%
%%%%%%%%%%%%%%%%%%%%%%%%%%%%%%%%%%%%%%%%%%%%%%%%%%%%%%%%%%%%%%%%%%%%%%%%%%%%%%
In particular, if $G$ is a vertex transitive graph,
$\chi_f(G)=\frac{|G|}{\alpha(G)}$ and hence
$\Sigma(G)<\frac{|G|^2}{\alpha(G)}$. Furthermore,
$e(G)=\frac{\Delta(G)|G|}{2}$. If
$\chi_f(G)\leq\frac{3}{4}\Delta(G)$, then $\chi_f(G)|G|<
\frac32(e(G)+1)$. Therefore, the bound in Theorem \ref{chif} is
better than the upper bound $\frac32(e(G)+1)$ (see
\cite{MR1000084}).

On the other hand, in view of Theorem \ref{hom} we have
$$\Sigma(G)\geq\frac{\omega(G)+1}{2}|G|$$ where $G$ is a vertex
transitive graph and $\omega(G)$ is the size of the largest
clique in it.

Also, It is a known result that the ratio
$\frac{\chi(G)}{\chi_f(G)}$ can be arbitrary large (see
\cite{MR1481157}). Let $\mathcal{G}=\{G_i\}_{i\in\mathbb{N}}$
such that $\frac{\chi(G_n)}{\chi_f(G_n)}\rightarrow\infty$. We
can assume that $G_n$ is critical  for all $n$ ($G$ is critical if
$\chi(G\setminus v)<\chi(G)$ for every $v\in V(G)$). Thus,
$e(G_n)\geq\frac{|G_n|(\chi(G_n)-1)}{2}$ and we also have
$\frac{\frac32(e(G_n)+1)}{\chi_f(G_n)|G_n|}\longrightarrow\infty$.
It means the bound in  Theorem \ref{chif} is better than the upper
bound $\frac32(e(G)+1)$  for the graphs in $\mathcal{G}$.

In Theorem \ref{chif} we used an upper bound  of
$\Sigma(KG(m,n))$, but we do not  know the exact value of
$\Sigma(KG(m,n))$.  The improvement of this upper bound  yields
an improvement in Theorem \ref{chif}.
\begin{preprb}
What is the exact value of $\Sigma(KG(m,n))$? Is it true that
$\Sigma(KG(m,n))={m\choose n
}(\frac{m+1}{n+1}-(\frac{n-1}{2n+2})\frac{{2n \choose n}}{{m
\choose n}})$?
\end{preprb}

\noindent{\bf Acknowledgment}

We would like to thank Hossein Hajiabolhassan and Moharam Nejad
Iradmusa for their useful comments.

%%%%%%%%%%%%%%%%%%%%%%%%%%%%%%%%%%%%%%%%%%%%%%%%%%%%%%%%%%%%%%%%%%%%%%%%%%%


\begin{thebibliography}{10}

\bibitem{MR1041612}
Paul Erd{\H{o}}s, Ewa Kubicka, and Allen~J. Schwenk.
\newblock Graphs that require many colors to achieve their chromatic sum.
\newblock In {\em Proceedings of the {T}wentieth {S}outheastern {C}onference on
  {C}ombinatorics, {G}raph {T}heory, and {C}omputing ({B}oca {R}aton, {FL},
  1989)}, volume~71, pages 17--28, 1990.

\bibitem{MR1746464}
H.~Hajiabolhassan, M.~L. Mehrabadi, and R.~Tusserkani.
\newblock Minimal coloring and strength of graphs.
\newblock {\em Discrete Math.}, 215(1-3):265--270, 2000.

\bibitem{tabular}
H.~Hajiabolhassan, M.~L. Mehrabadi, and R.~Tusserkani.
\newblock Tabular graphs and chromatic sum.
\newblock {\em Discrete Math.}, 304(1-3):11--22, 2005.

\bibitem{MR2089014}
Pavol Hell and Jaroslav Ne{\v{s}}et{\v{r}}il.
\newblock {\em Graphs and homomorphisms}, volume~28 of {\em Oxford Lecture
  Series in Mathematics and its Applications}.
\newblock Oxford University Press, Oxford, 2004.

\bibitem{MR1722229}
Tao Jiang and Douglas~B. West.
\newblock Coloring of trees with minimum sum of colors.
\newblock {\em J. Graph Theory}, 32(4):354--358, 1999.

\bibitem{MR2085077}
Ewa Kubicka.
\newblock The chromatic sum of a graph: history and recent developments.
\newblock {\em Int. J. Math. Math. Sci.}, (29-32):1563--1573, 2004.

\bibitem{kobik}
Ewa Kubicka and Allen~J. Schwenk.
\newblock An introduction to chromatic sums.
\newblock {\em Proc. ACM Computer Science Conference},
  Louisville(Kentucky):39--45, 1989.

\bibitem{MR514625}
L.~Lov{\'a}sz.
\newblock Kneser's conjecture, chromatic number, and homotopy.
\newblock {\em J. Combin. Theory Ser. A}, 25(3):319--324, 1978.

\bibitem{MR1481157}
Edward~R. Scheinerman and Daniel~H. Ullman.
\newblock {\em Fractional graph theory}.
\newblock Wiley-Interscience Series in Discrete Mathematics and Optimization.
  John Wiley \& Sons Inc., New York, 1997.
\newblock A rational approach to the theory of graphs, With a foreword by
  Claude Berge, A Wiley-Interscience Publication.

\bibitem{MR1000084}
Carsten Thomassen, Paul Erd{\H{o}}s, Yousef Alavi, Paresh~J.
Malde, and
  Allen~J. Schwenk.
\newblock Tight bounds on the chromatic sum of a connected graph.
\newblock {\em J. Graph Theory}, 13(3):353--357, 1989.

\end{thebibliography}
\end{document}